\documentclass[10.5pt]{amsart}
\usepackage{amsmath}
\usepackage[applemac]{inputenc}
\usepackage[frenchb]{babel}
\usepackage{amscd}
\usepackage{xspace}
\usepackage{verbatim}

\newcommand{\fnz}{\footnotesize}

\newcommand{\be}{\begin{equation}}
\newcommand{\bel}[1]{\begin{equation}\label{#1}}
\newcommand{\ee}{\end{equation}}
\newtheorem{subn}{\name}

\newcommand{\bsn}[1]{\def\name{#1}\begin{subn}}
\newcommand{\esn}{\end{subn}}
\newtheorem{sub}{\name}
\newcommand{\bs}{\begin{sub}}
\newcommand{\es}{\end{sub}}
\newcommand{\bsl}[1]{\begin{sub}\label{#1}}
\newcommand{\bth}[1]{\def\name{Theorem}\begin{sub}\label{t:#1}}
\newcommand{\blemma}[1]{\def\name{Lemma}\begin{sub}\label{l:#1}}
\newcommand{\bcor}[1]{\def\name{Corollary}\begin{sub}\label{c:#1}}
\newcommand{\bdef}[1]{\def\name{Definition}\begin{sub}\label{d:#1}}
\newcommand{\bprop}[1]{\def\name{Proposition}\begin{sub}\label{p:#1}}

\newcommand{\BA}{\begin{array}}
\newcommand{\EA}{\end{array}}
\newcommand{\BAN}{\renewcommand{\arraystretch}{1.2}
\setlength{\arraycolsep}{2pt}\begin{array}}
\newcommand{\BAV}[2]{\renewcommand{\arraystretch}{#1}
\setlength{\arraycolsep}{#2}\begin{array}}
\newcommand{\BSA}{\begin{subarray}}
\newcommand{\ESA}{\end{subarray}}
\newcommand{\BAL}{\begin{aligned}}
\newcommand{\EAL}{\end{aligned}}
\newcommand{\BALG}{\begin{alignat}}
\newcommand{\EALG}{\end{alignat}}
\newcommand{\BALGN}{\begin{alignat*}}
\newcommand{\EALGN}{\end{alignat*}}

\newcommand{\abs}[1]{\left |#1\right |}
\newcommand{\norm}[1]{\left \|#1\right \|}

\newcommand{\opname}[1]{\mbox{\rm #1}\,}

\newcommand{\dist}{\opname{dist}}
\newcommand{\myfrac}[2]{{\displaystyle \frac{#1}{#2} }}
\newcommand{\myint}[2]{{\displaystyle \int_{#1}^{#2}}}

\newcommand{\prt}{\partial}

\def\ga{\alpha}            
       \def\gd{\delta}      
\def\gth{\theta}                         
           
      \def\gk{\kappa}      \def\gl{\lambda}
\def\gm{\mu}        \def\gn{\nu}         \def\gp{\pi}
    \def\gr{\rho}        
\def\gs{\sigma}       
      \def\gw{\omega}
                \def\gz{\zeta}
     \def\Gd{\Delta}      

    \def\Gs{\Sigma}      
\def\Gw{\Omega}              

\def\CS{{\mathcal S}}      
\def\CR{{\mathcal R}}      
      
   \def\CE{{\mathcal E}}   \def\CF{{\mathcal F}}
      
\def\CJ{{\mathcal J}}      \def\CL{{\mathcal L}}

\def\BBG {\mathbb G}       
   \def\BBK {\mathbb K}    
       
   \def\BBR {\mathbb R}

\def\GTM {\mathfrak M}

\begin{document}%

\noindent Partial Differential Equations\smallskip

\title[ Semilinear equations with Hardy potentials]{Measure boundary value problem for semilinear elliptic equations with critical Hardy potentials}
\maketitle
\noindent {Konstantinos T. Gkikas
\footnote{\noindent Centro de Modelamiento Matem\`atico,
Universidad de Chile, Santiago de Chile, Chile. E-mail: kugkikas@gmail.com},
 Laurent V\'{e}ron
 \footnote{\noindent 
Laboratoire de Math\'{e}matiques et Physique Th\'{e}orique, CNRS UMR 7350,
Facult\'{e} des Sciences, 37200 Tours France. E-mail: veronl@univ-tours.fr}}

\date{}
\begin{abstract} Let $\Omega\subset\BBR^N$ be a bounded $C^2$ domain and $\CL_\gk=-\Gd-\frac{\gk}{d^2}$ the Hardy operator where $d=\dist (.,\prt\Gw)$ and $0<\gk\leq\frac{1}{4}$. Let $\ga_{\pm}=1\pm\sqrt{1-4\gk}$ be the two Hardy exponents,  $\gl_\gk$ the first eigenvalue of $\CL_\gk$ with corresponding positive eigenfunction $\phi_\gk$. If $g$ is a continuous nondecreasing function satisfying $\int_1^\infty(g(s)+|g(-s)|)s^{-2\frac{2N-2+\ga_+}{2N-4+\ga_+}}ds<\infty$, 
then for any Radon measures $\gn\in \GTM_{\phi_\gk}(\Gw)$ and $\gm\in \GTM(\prt\Gw)$ there exists a unique weak solution
 to problem $P_{\gn,\gm}$: $\CL_\gk u+g(u)=\gn$ in $\Gw$, $u=\gm$ on $\prt\Gw$. If $g(r)=|r|^{q-1}u$ ($q>1$) we prove that, in the subcritical range of $q$, a necessary and sufficient condition for solving $P_{0,\gm}$ with $\gm>0$ is that $\gm$ is absolutely continuous with respect to the capacity associated to the Besov space $B^{2-\frac{2+\ga_+}{2q'},q'}(\BBR^{N-1})$. We also characterize the boundary removable sets in terms of this capacity. In the subcritical range of $q$ we classify the isolated singularities of positive solutions.
  \end{abstract}
\begin{center}{\bf Problèmes aux limites avec données mesures pour des équations semi linéaires elliptiques avec des potentiels de Hardy critiques }\end{center}
\smallskip
\begin{quotation}
 {\fnz{\sc R\'esum\'e.} Soient $\Omega\subset\BBR^N$ un domaine de classe $C^2$ et $\CL_\gk=-\Gd-\frac{\gk}{d^2}$ l'opérateur de Hardy où $d=\dist (.,\prt\Gw)$ et $0<\gk\leq\frac{1}{4}$. Soient $\ga_{\pm}=1\pm\sqrt{1-4\gk}$ les deux exposants de Hardy,  $\gl_\gk$ première valeur propre de $\CL_\gk$ et $\phi_\gk$ la fonction propre positive correspondante.  Si $g$ est une fonction continue croissante vérifiant $\int_1^\infty(g(s)+|g(-s)|)s^{-2\frac{2N-2+\ga_+}{2N-4+\ga_+}}ds<\infty$, alors pour toutes mesures de Radon $\gn\in \GTM_{\phi_\gk}(\Gw)$ et $\gm\in \GTM(\prt\Gw)$ il existe une unique solution faible
au problème $P_{\gn,\gm}$: $\CL_\gk u+g(u)=\gn$ dans $\Gw$, $u=\gm$ sur $\prt\Gw$. Si $g(r)=|r|^{q-1}u$ ($q>1$) nous démontrons qu'une condition nécessaire et suffisante pour résoudre $P_{0,\gm}$ avec $\gm>0$ est que $\gm$ soit absolument continue par rapport à la capacité associée à l'espace de Besov $B^{2-\frac{2+\ga_+}{2q'},q'}(\BBR^{N-1})$. Nous caractérisons les ensembles éliminables pour les valeurs sur critiques de $q$. Dans le cas sous -critique nos donnons une classifications des singularités isolées au bord des solutions positives.
}
\end{quotation}
 \bigskip
 {\bf Version fran\c{c}aise abr\'eg\'ee.}
Soit $\Gw$ un domiane  de $\BBR^N$ de classe $C^2$. On désigne par $d(x)$ la distance de $x$ à $\prt\Gw$ et on définit  l'opérateur de Hardy dans $\Gw$ par
\bel{F1}
\CL_\gk u=-\Gd u-\frac{\gk}{d^2}u
\ee
où $0<\gk\leq\frac{1}{4}$ et ses exposants caractéristiques  
\bel{F2}\ga_+=1+\sqrt{1-4\gk}\quad\ga_-=1-\sqrt{1-4\gk}.\ee
On supposera $\Gw$ convexe si $\gk=\frac{1}{4}$. Il est bien connu que sous ces conditions $\CL_\gk$ possède une première valeur propre $\gl_\gk>0$ définie par
\bel{F3}\displaystyle
\gl_\gk:=\inf_{u\in H^1_0(\Gw)\setminus\{0\}}\myfrac{\myint{\Gw}{}|\nabla u|^2dx}{\myint{\Gw}{}d^{-2}u^2dx}.
\ee
La  première fonction propre positive associée $\phi_\gk$ n'appartient à $H^1_0(\Gw)$ que si $0<\gk<\frac{1}{4}$, et dans tous les cas elle vérifie $\phi_\gk(x)\sim (d(x))^{\ga_+}$ au voisinage de $\prt\Gw$. On dénote par $G_\gk$ et $K_\gk$ les noyaux de Green et de Poisson de $\CL_\gk$ dans $\Gw$ et par $\gw^{x_0}$ la mesure harmonique dans $\Gw$ ($x_0\in\Gw$).
Si $g$ est une fonction continue et croissante sur $\BBR$ telle que $g(0)\geq 0$, nous étudions tout d'abord le problème $(P_{\gn,\gm})$ suivant
\bel {F4}\BA {lll}
\CL_\gk u+g(u)=\gn\qquad&\text{in }\Gw\\
\phantom{\CL_\gk +g(u)}
u=\gm\qquad&\text{in }\prt\Gw,
\EA\ee
où $\gn$, $\gm$ sont des mesures de Radon. \medskip

\noindent{\bf Théorème 1}. {\it Supposons que $g$ vérifie
\bel{F5}\BA {lll}
\myint{1}{\infty}\left(g(s)+|g(-s)|\right) s^{-2\frac{N-1+\frac{\ga_+}{2}}{N-2+\frac{\ga_+}{2}}}ds<\infty;
\EA\ee
alors pour toutes mesures de Radon $\gn$ et $\gm$ dans $\Gw$ et $\prt\Gw$ respectivement, $\gn$ vérifiant en outre 
$\int_\Gw\phi_\gk d|\gn|<\infty$, il existe une unique fonction $u=u_{\gn,\gm}\in L^1_{\phi_\gk}(\Gw)$ telle que $g{\footnotesize {\circ}}u\in L^1_{\phi_\gk}(\Gw)$ vérifiant 
\bel{F6}\BA {lll}
\myint{\Gw}{}\left(u\CL_\gk\gz+\gz g{\circ}u\right)dx=\myint{\Gw}{}G_\gk(x,y) d\gn(y)+\myint{\prt\Gw}{}K_\gk(x,y) d\gm(y)
\EA\ee
pour toute $\gz\in {\bf X}_\gk(\Gw)$ où 
\bel{F7}\BA {lll}
{\bf X}_\gk(\Gw)=\{\gz\in H^1_{loc}(\Gw):(\phi_\gk)^{-1}\in H^1_{0}(\Gw,\phi_\gk dx), (\phi_\gk)^{-1}\CL_\gk\in L^\infty(\Gw)\}.
\EA\ee
En outre l'application $(\gn,\gm)\mapsto u_{\gn,\gm}$ de $\mathfrak M_{\phi_\gk}(\Gw)\times \mathfrak M(\prt\Gw)$ dans $L^1_{\phi_\gk}(\Gw)$ est croissante et stable pour la convergence faible des mesures.}\medskip

La démonstration utilise des estimations des noyaux de Green et de Poisson obtenus à partir des propriétés de la mesure harmonique. Dans le cas où $g{\circ}u=\abs u^{q-1}u$ l'inégalité (\ref{F5}) est vérifiée si $0<q<q_c:=\frac{2N+\ga_+}{2N+\ga_+-4}$
Dans le cas $q>1$ nous dénotons par $C^{\BBR^{N-1}}_{2-\frac{2+\ga_+}{2q'},q'}$ la capacité associée à  l'espace de Besov $B^{2-\frac{2+\ga_+}{2q'},q'}(\BBR^{N-1})$ et nous avons le résultat suivant:\medskip

\noindent{\bf Théorème 2}. {\it Soit $q\geq q_c$ et $\gn\in \mathfrak M_+(\prt\Gw)$. Alors le problème
\bel {F8}\BA {lll}
\CL_\gk u+\abs u^{q-1}u=0\qquad&\text{in }\Gw\\
\phantom{\CL_\gk +\abs u^{q-1}u}
u=\gm\qquad&\text{in }\prt\Gw
\EA\ee
admet une unique solution $u:=u_\gm$ si et seulement si pour tout borélien $E\subset\prt\Gw$, 
\bel {F9}\BA {lll}
C^{\BBR^{N-1}}_{2-\frac{2+\ga_+}{2q'},q'}(E)=0\Longrightarrow \gm (E)=0.
\EA\ee
}

Nous caractérisons aussi les sous ensembles du bord éliminables pour 
\bel{F8'}
\CL_\gk u+\abs u^{q-1}u=0\qquad \text{in }\Gw.
\ee
Pour cela nous posons
\bel {F11}
W(x)=\left\{\BA {lll} &(d(x))^{\frac{\ga_-}{2}}\qquad&\text{if }0<\gk<\frac{1}{4}\\[2mm]
&\sqrt{d(x)}\ln|d(x)|\qquad&\text{if }\gk=\frac{1}{4}.
\EA\right.\ee
\medskip

\noindent{\bf Théorème 3}. {\it Soit $q>1$ et $K\subset\prt\Gw$ un sous-ensemble compact. Toute solution $u\in C(\overline\Gw\setminus \{K\})$ de (\ref{F8'}) qui vérifie
\bel {F12}\BA {lll}
\lim_{x\to y}\myfrac{u(x)}{W(x)}=0\qquad \forall y\in\prt\Gw\setminus\{K\},
\EA\ee
est identiquement nulle dans $\Gw$ si et seulement si $C^{\BBR^{N-1}}_{2-\frac{2+\ga_+}{2q'},q'}(K)=0$.
}\medskip

Nous montrons que si $q>1$, toute solution positive de (\ref{F8'}) dans $\Gw$ admet une trace au bord représentée par une mesure de Borel régulière. En supposant que $0\in\prt\Gw$ et $1<q<q_c$, nous étudions aussi le comportement au voisinage de $0$ des solutions positives de (\ref{F8'}) qui vérifient (\ref{F12}) avec $K=\{0\}$.
\begin {center}
---------------------------------------------------------------------------------------------
\end {center} 

Let $\Gw$ be a bounded $C^2$ domain in $\BBR^N$ and $d(x)=\dist (x,\Gw)$. We define the Hardy operator $\CL_\gk$ in 
$\Gw$ by (\ref{F1}) with $0<\gk\leq \frac{1}{4}$ and the characteristic exponents by (\ref{F2}). We assume that $\Gw$ is convex if 
$\gk=\frac{1}{4}$. It is well known that $\CL_\gk$ possesses a first eigenvalue $\gl_\gk>0$ defined by (\ref{F3}) and that the first positive eigenfunction $\phi_\gk>0$ may or may not belong to $H^1_0(\Gw)$ according $0<\gk<\frac{1}{4}$ or $\gk=\frac{1}{4}$,
and it satisfies $\phi_\gk(x)\sim (d(x))^{\frac{\ga_+}{2}}$, $|\nabla\phi_\gk(x)|\sim (d(x))^{\frac{\ga_+}{2}-1}$ as $d(x)\to 0$. The Green and the Poisson kernels are denoted by $G_\gk(x,y)$ and $K_\gk(x,y)$, and they satisfy
\bel {F13}
G_\gk(x,y)\sim  \min\left\{\myfrac{1}{|x-y|^{N-2}},\frac{(d(x))^{\myfrac{\ga_+}{2}}(d(y))^{\frac{\ga_+}{2}}}{|x-y|^{N-2+\ga_+}}\right\}\quad \forall (x,y)\in \Gw\times\Gw,\,x\neq y,
\ee
\bel {F14}
K_\gk(x,y)\sim  \myfrac{(d(x))^{\frac{\ga_+}{2}}}{|x-y|^{N-2+\ga_+}}\quad \forall (x,y)\in \Gw\times\prt\Gw.
\ee
The corresponding Green and Poisson operators are denoted by $\BBG_\gk[.]$ and $\BBK_\gk[.]$. We first consider the boundary value problem (\ref{F4}) where $g$ is a continuous nondecreasing function such that $g(0)\geq 0$ and $\gn$ and $\gm$ are Radon measures in $\Gw$ and $\prt\Gw$ respectively. We say that $g$ is a subcritical nonlinearity if it satisfies (\ref{F5}).
 \medskip

\bth{1}Assume that $g$ is a subcritical nonlinearity.
Then for all $(\gn,\gm)\in \mathfrak M_{\phi_\gk}(\Gw)\times \mathfrak M(\prt\Gw)$ there exists a unique function $u=u_{\gn,\gm}\in L^1_{\phi_\gk}(\Gw)$ such that $g\footnotesize {\circ}u\in L^1_{\phi_\gk}(\Gw)$  verifying (\ref{F6}) for all $\gz$ in the space of test functions ${\bf X}_\gk(\Gw)$ defined by (\ref{F7}).
Furthermore the mapping $(\gn,\gm)\mapsto u_{\gn,\gm}$ from $\mathfrak M_{\phi_\gk}(\Gw)\times \mathfrak M(\prt\Gw)$ into $L^1_{\phi_\gk}(\Gw)$ is nondecreasing and stable for the weak convergence of measures.\es

When $g(u)=|u|^{q-1}u$ with $q>0$, the inequality (\ref{F6}) means  
 \bel {F15}
0<q<q_c:=\frac{2N+\ga_+}{2N+\ga_+-4}.
\ee
When $q\geq q_c$ not all the measures $\gm$ are eligible for solving (\ref{F8}). We denote by $C^{\BBR^{N-1}}_{2-\frac{2+\ga_+}{2q'},q'}$ the capacity associated to the Besov space $W^{2-\frac{2+\ga_+}{2q'},q'}(\BBR^{N-1})$.

\bth{2}Let $q>1$ and $\gn\in \mathfrak M_+(\prt\Gw)$. Then problem (\ref{F8})
admits a solution if and only if $\gm$ is absolutely continuous with respect to  $C^{\BBR^{N-1}}_{2-\frac{2+\ga_+}{2q'},q'}$, i.e. for any Borel set $E\subset\prt\Gw$, implication (\ref{F9}) holds.
\es

We also characterize the boundary removable sets for (\ref{F8'}).
\bth{3} Let $q>1$ and $K\subset\prt\Gw$ is compact. Any $u\in C(\overline\Gw\setminus \{K\})$ solution of 
(\ref{F8'}) which verifies (\ref{F12})
is identically zero in $\Gw$ if and only if $C^{\BBR^{N-1}}_{2-\frac{2+\ga_+}{2q'},q'}(K)=0$.
\es

When $1<q<q_c$ only the empty set has zero capacity. There exist  singular solutions of (\ref{F8'}) with an isolated singularity on the boundary either solutions $u_{k\gd_a}$ of (\ref{F8}) with $\gm=k\gd_a$ for $k>0$ and $a\in\prt\Gw$ or solutions 
$u_{a}=\lim_{k\to\infty}u_{k\gd_a}$. This very singular solution is described by considering the following problem on 
the half sphere $S^{N-1}_+=\{x=(x_1,...,x_N)\in\BBR^N:\abs x=1,\,x_N=1\}$
 \bel {F16}\BA {ll}
 -\Gd'\gw-\ell_{N,q,\gk}\gw-\frac{\gk}{({\bf e}_N.\gs)^2}\gw+|\gw|^{q-1}\gw=0\qquad&\text{in }\,S^{N-1}_+\\
 \phantom{-\Gd'\gw-\ell_{N,q,\gk}\gw-\frac{\gk}{({\bf e}_N.\gs)^2}\gw+|\gw|^{q-1}}
 \gw=0\qquad&\text{in }\,\prt S^{N-1}_+
\EA\ee 
where $\Gd'$ is the Laplace-Beltrami operator on $S^{N-1}$, $({\bf e}_1,...,{\bf e}_N)$ is the canonic basis in $\BBR^N$, $\gs=\frac{x}{|x|}$ and 
$$\ell_{N,q}=\left(\frac{2}{q-1}\right)\left(\frac{2q}{q-1}-N\right).
$$

The spherical Hardy  operator $\gw\mapsto \CL_\gk':= -\Gd'\gw-\frac{\gk}{({\bf e}_N.\gs)^2}\gw$ on $S^{N-1}_+$ admits a first eigenvalue $\gm_\gk$ defined by
 \bel {F17}\BA {ll}
\gm_{\gk,1}=\displaystyle\inf_{\psi\in H^1_0(S^{N-1}_+)\setminus\{0\}}
\myfrac{\myint{S^{N-1}_+}{}\left(|\nabla' \psi|^2-\gk({\bf e}_N.\gs)^{-2}\gw^2\right)dS}{\myint{\Gw}{}({\bf e}_N.\gs)^{-2}\psi^2dS}.
\EA\ee
We prove that 
$\gm_{\gk,1}=\frac{\ga_+}{2}\left(N+\frac{\ga_+}{2}-2\right)$
with corresponding positive eigenfunction $\gr_{\gk}=({\bf e}_N.\gs)^{\frac{\ga_+}{2}}$. There exists a second eigenvalue $\gm_{\gk,2}=\gm_{\gk,1}+N+\ga_+-1$ with $N-1$ independent eigenfunctions 
$\gr_{\gk,j}=({\bf e}_N.\gs)^{\frac{\ga_+}{2}}{\bf e}_j.\gs$ for $j=1,...,N-1$.  We denote by $\CE_{\gk}$ the set of functions $\gw$ such that 
$\gr_\gk^{-1}\gw\in L^{q+1}_{\gr_\gk^{q+1}}(S^{N-1}_+)\cap H^{1}_0(S^{N-1}_+,\gr_\gk^2dS)$ which satisfy (\ref{F16}), and by $\CE^+_{\gk}$ the set of positive solutions.
\bth{NJEV} I- If $q\geq q_c$, $\CE_{\gk}=\{\emptyset\}$.\\
II- If $1<q< q_c$, $\CE^+_{\gk}=\{0,\gw\gk\}$ where $\gw_\gk$ is the unique positive solution of (\ref{F16}).\\
III- If $q_e\leq q<q_c$, $\CE_{\gk}=\{0,\gw_\gk,-\gw_\gk\}$ where $q_e:=\myfrac{2N+2+\ga_+}{2N-2+\ga_+}$.
\es

This result allows to describe the isolated  boundary singularities of positive solutions of (\ref{F8'}). We assume that $0\in\prt\Gw$ and the outward normal unit vector to $\prt\Gw$ at $0$ is ${\bf e}_N$.

\bth{Class} Assume ,  $1<q<_c$ and $u\in C(\overline\Gw\setminus\{0\})$ is a positive solution of (\ref{F8'}) which verifies(\ref{F12}) with $K=\{0\}$. Then\\
(i) either there exists $k\geq 0$ such that $u=u_{k\gd_0}$ and $\lim_{|x|\to 0}|x|^{N+\frac{\ga_+}{2}-2}u(x)=c_Nk({\bf e}_N.\frac{x}{|x|})^{\frac{\ga_+}{2}}$,\\
(ii) or $\lim_{|x|\to 0}|x|^{\frac{2}{q-1}}u(x)=\gw_\gk(\frac{x}{|x|})$.\\
The above two convergence hold  locally uniformly on $S^{N-1}_+$.
\es

We can also define a boundary trace of any positive solution $u$ of (\ref{F8'}). For $\gd>0$ small enough, we 
denote by $\gw_{\Gw'_{\gd}}^{x_0}$ the harmonic measure relative to the operator $\CL_\gk$ in $\Gw'_{\gd}=\{x\in\Gw:d(x)>\gd\}$
where $x_0\in \Gw$ (with $d(x_0)\geq\gd_1>\gd$) and set $\Gs_\gd=\prt\Gw'_{\delta}$.

\bth{Tr} Assume  $q>1$ and $u\in C(\overline\Gw\setminus\{0\})$ is a positive solution of (\ref{F8'})  in $\Gw$. Then for any $y\in\prt\Gw$, the following dichotomy occurs:\\
(i) Either there exists an open subset $U\subset\BBR^N$ containing $y$ and a positive Radon measure $\gl_U$ on $\prt\Gw\cap U$ such that 
 \bel {T1}
 \lim_{\gd\to 0}\myint{\Gs_\gd\cap U}{}Z(x)u(x)d\gw_{\Gw'_{\gd}}^{x_0}=\myint{\prt\Gw\cap U}{}Zd\gl_U\quad\forall Z\in C_0(U).
\ee
(ii) Or for any open subset $U\subset\BBR^N$ containing $y$, there holds
 \bel {T2}
 \lim_{\gd\to 0}\myint{\Gs_\gd\cap U}{}u(x)d\gw_{\Gw'_{\gd}}^{x_0}=\infty.
\ee
\es

The set $\CR_u$ of $x_0$ such that (i) holds is relatively open in $\prt\Gw$ and it carries a positive Radon measure $\gm_u$ such that 
(\ref{T1}) occurs with $U$ replaced by $\CR_u$ and $\gl_U$ by $\gm_u$; its complement  $\CS_u$ in $\prt\Gw$ has the property that (\ref{T2}) occurs for any open subset $U$ such that $U\cap \CS_u\neq\{\emptyset\}$.

\noindent{\it Abridged proof of Theorem 1.}  Let $(\gn,\gm)\in \mathfrak M_{\phi_\gk}(\Gw)\times  \mathfrak M(\prt\Gw)$. For $\gl>0$ we set 
 \bel {F18}E_\gl(\gn)=\{x\in\Gw:\BBG_\gk[|\gn|](x)>\gl\},\;\CE_\gl(\gn)=\myint{E_\gl(\gn)}{}\phi_\gk dx,
\ee
and
 \bel {F10}F_\gl(\gn)=\{x\in\Gw:\BBK_\gk[|\gm|](x)>\gl\},\;\CF_\gl(\gm)=\myint{E_\gl(\gn)}{} dx,
\ee
and prove
 \bel {F20}
 \CE_\gl(\gn)+\CF_\gl(\gm)\leq c\left(\myfrac{\norm\gn_{\mathfrak M_{\phi_\gk}(\Gw)}+
 \norm\gm_{\mathfrak M(\prt\Gw)}}{\gl}\right)^{\frac{2N+\ga_+}{2N+\ga_+-4}}.
 \ee
If $g$ satisfies (\ref{F5}) and $\{(\gn_n,\gm_n)\}$ is a sequence of smooth functions which converges in the weak-star topology of measures to 
$(\gn,\gm)$, then the corresponding solutions $\{u_{\gn_n,\gm_n}\}$ of problem $P_{\gn_n,\gm_n}$ defined in  ((\ref{F4})) converges to some $u$ and $\{g\circ u_{\gn_n,\gm_n}\}$ converges to $g\circ u$ in $L^1_{\phi_\gk}$ by Vitali convergence theorem. This implies that $u=u_{\gn,\gm}$. Uniqueness holds by adapting Brezis estimates and using monotonicity.\medskip

\noindent {\it Abridged proof of Theorem 2.} Using estimate (\ref{F14}) and the harmonic lifting in Besov spaces introduced in \cite[Sect. 3] {MV2} we prove that for any $\gm\in \mathfrak M(\prt\Gw)$ there holds
 \bel {F21}
\frac{1}{c}\norm{\gm}^q_{B^{-2+\frac{2+\ga_+}{2q'},q}}\leq \myint{\Gw}{}(\BBK[|\gm|])^q\phi_\gk
dx\leq  c\norm{\gm}^q_{B^{-2+\frac{2+\ga_+}{2q'},q}}\ee
for some $c=c(\Gw,\gk,q)>0$. This implies that we can solve (\ref{F8}) with such a Radon measure.  If $\gm\in B^{-2+\frac{2+\ga_+}{2q'},q}(\prt\Gw)\cap \mathfrak M_+(\prt\Gw)$, it is absolutely continuous with respect to the capacity $C^{\BBR^{N-1}}_{2-\frac{2+\ga_+}{2q'},q}$. Finally, if $\gm\in \mathfrak M_+(\prt\Gw)$ is absolutely continuous with respect to the capacity $C^{\BBR^{N-1}}_{2-\frac{2+\ga_+}{2q'},q}$, there exists an increasing sequence $\{\gm_n\}\subset B^{-2+\frac{2+\ga_+}{2q'},q}(\prt\Gw)\cap \mathfrak M_+(\prt\Gw)$ which converges to $\gm$. This implies that  $u_{\gm_n}$ converges to $u_{\gm} $ in $L^q_{\phi_\gk}(\Gw)$.\\
Conversely, if  $\gm\in \mathfrak M_+(\prt\Gw)$ is such that there exists a solution $u_{\gm}$ to (\ref{F8}), we use a variant of the {\it optimal lifting} $R[.]$ defined in \cite[Sect. 1]{MV-JMPA} to prove that for any $\eta\in C^2(\prt\Gw)$ such that $0\leq \eta\leq 1$ there holds
 \bel {F22}
\myint{\prt\Gw}{}\eta d\gm\leq c\myint{\Gw}{}u^q\gz dx+c\left(\myint{\Gw}{}u^q\gz dx\right)^{\frac{1}{q}}
\left(\myint{\Gw}{}\phi_\gk dx+\norm{\eta}^{q'}_{B^{2-\frac{2+\ga_+}{2q'},q"}}\right)^{\frac{1}{q'}}.
\ee
Here $\gz=\phi_\gk (R[\eta])^{q'}$ and $R:C^2(\prt\Gw)\mapsto C^2(\overline\Gw)$ is a linear mapping which satisfies 
$0\leq \eta\leq 1\Longrightarrow 0\leq R[\eta]\leq 1$ and $R[\eta]\lfloor_{\prt\Gw}=\eta.$ If $K\subset\prt\Gw$ is a compact set
with zero $C^{\BBR^{N-1}}_{2-\frac{2+\ga_+}{2q'},q}$-capacity, there exists a sequence $\{\eta_n\}\subset C^2(\prt\Gw)$ such that $0\leq\eta_n\leq 1$, $\eta_n=1 $ on $K$ and $\norm{\eta_n}^{q'}_{B^{2-\frac{2+\ga_+}{2q'},q"}}\to 0$. This implies 
$\phi_\gk (R[\eta_n])^{q'}\to 0$ and finally $\gm (K)=0$.\medskip

\noindent {\it Abridged proof of Theorem 3.} If $K\subset\prt\Gw$ is compact with $C^{\BBR^{N-1}}_{2-\frac{2+\ga_+}{2q'},q}(K)>0$, its capacitary measure $\gm_K$ belongs to $B^{-2+\frac{2+\ga_+}{2q'},q}(\prt\Gw)\cap \mathfrak M_+(\prt\Gw)$ . Thus $u_{\gm_K}$ exists and $K$ is not removable. Conversely by using again optimal lifting, and test functions of the form
$\phi_\gk (R[1-\eta])^{2q'}$ where $0\leq\eta\leq 1$ and $\eta=1 $ in a neighborhood of $K$, we prove first that that $u\in L^q_{\phi_\gk}(\Gw)$ and finally that $u=0$.\medskip

\noindent {\it Abridged proof of Theorems 4-5.} Existence is obtained in minimizing $\CJ_\gk$ defined over 
$L^{q+1}_{\gr_\gk^{q+1}}(S^{N-1}_+)\cap H^{1}_0(S^{N-1}_+,\gr_\gk^2dS)$ by
 \bel {F23}
\CJ_\gk(w):=\myint{S^{N-1}_+}{}\left(|\nabla'w|^2-(\ell_{N,q}-\gm_{\gk,1})w^2+\frac{2}{q+1}\gr_\gk^{q-1}|w|^{q+1}\right)\gr_\gk^2 dS.
\ee
A non-trivial minimizer exists if $\ell_{N,q}>\gm_{\gk,1}$ (defined by (\ref{F17})), i.e. $1<q<q_c$, and $\gw=\gr_\gk w$ satisfies (\ref{F16}).
Nonexistence in standard since $\gm_{\gk,1}<\ell_{N,q}$ if and only if $1<q<q_c$. For uniqueness we assume that $\gw_j$ ($j=1,2$) are positive solutions of (\ref{F16}) and we set
$w_j=\frac{\gw_j}{\gr_\gk}$. Then
$$-{\rm div'}.(\gr_\gk^2\nabla' w_j)+(\gm_{\gk,1}-\ell_{N,q})\gr_\gk^2w_j+\gr_\gk^{q+1}w_j^q=0\qquad\text{on }S^{N-1}_+
$$
Near $\prt S^{N-1}_+$ we have $w_j\sim \gr_\gk^{\frac{\ga_+}{2}}$ and $|\nabla'w_j|\sim \gr_\gk^{\frac{\ga_+}{2}-1}$. Then integration by parts is justified and

$$\myint{S^{N-1}_+}{}\left(\left(\myfrac{\nabla' w_1}{w_1}-\myfrac{\nabla' w_2}{w_2}\right).\nabla' (w^2_1-w^2_2)
+\gr_\gk^{q-1}(w_1^{q-1}-w_2^{q-1})(w^2_1-w^2_2)\right)\gr_\gk^2dS=0.
$$
The two terms of the integral are nonnegative, thus $w_1=w_2$. For statement III we first prove, by the method used in \cite[Th 3.1]{Ve89}, that any solution $\gw$ depends only on the azimuthal angle $\gth\in [0,\frac{\gp}{2}]$. Then we show that the corresponding ODE verified by $\gw$ admits only the three mentioned solutions. For Theorem 5, we first construct a barrier function as in \cite[Appendix]{MV-CPAM} which yields to  the following estimate
 \bel {F24}
u(x) \leq c|x|^{-\frac{2}{q-1}+\frac{\ga_+}{2}}(d(x))^{\frac{\ga_+}{2}}\qquad\forall x\in\Gw.
\ee
With this estimate we adapt the scaling method developed in \cite[Sect. 3.3]{NPV} to obtain the classification result. 


\begin{thebibliography}{99}

\bibitem{1} C.Bandle, V.Moroz and W.Reichel, \emph{Boundary blow up type sub-solutions to semilinear elliptic
equations with Hardy potential}, J. London Math. Soc. 77, 503-523 (2008).

\bibitem{BFT1} G. Barbatis, S. Filippas and A. Tertikas, \emph{A unified approach to improved $L^p$ Hardy inequalities with best constants.} Trans. Amer. Math. Soc. 356, 2169–2196 (2003).
\bibitem{BM} H. Brezis and M. Marcus, \emph{Hardy's inequalities revisited } Ann. Sc. Norm. Super. Pisa Cl. Sci. 25 (1997), no. 4, 217-237.  
\bibitem{dd} J. Davila, L. Dupaigne,{ \it Hardy-type inequalities,} J. Eur. Math. Soc. (JEMS) 6 (3) (2004) 335-365.
\bibitem{F.M.T2} S. Filippas, L. Moschini. and A. Tertikas, \emph{Sharp two-sided heat kernel estimates for critical Schrodinger operators on bounded domains.} Comm. Math. Phys. 273 (2007), 237-281. 
\bibitem{MV-JMPA} M. Marcus, L. V\'{e}ron, {\it Removable singularities and boundary trace.} J. Math. Pures Appl. 80, 879-900  (2001). 
\bibitem{MV-CPAM} M. Marcus, L. V\'{e}ron, {\it  The boundary trace and generalized boundary value problem for semilinear elliptic equations with coercive absorption.} Comm. Pure Appl. Math. 56 (2003), 689-731. 
\bibitem{MV2} M. Marcus, L. V\'{e}ron, {\it Boundary trace of positive solutions of supercritical semilinear elliptic equations in dihedral domains.} Ann. Sc. Norm. Super. Pisa Cl. Sci., to appear. arXiv:1309.7778.
\bibitem{NPV}Nguyen Phuoc T., L. V\'{e}ron, {\it Boundary singularities of solutions to elliptic viscous Hamilton–Jacobi equations.} J. Funct. Anal. 263 (2012) 1487–1538.
\bibitem{Ve89} L. V\'{e}ron, {\it Geometric invariance of singular solutions of some nonlinear partial differential equations.} Indiana Univ. Math. J. 38 (1989), 75-100.
\end{thebibliography}
\end {document}